\documentclass[12pt,a4paper,fleqn]{article}
\usepackage{a4wide,amsfonts,amsmath,latexsym,amssymb,euscript,graphicx,units,mathrsfs}

\usepackage[english]{babel}
\usepackage{graphicx}
\usepackage{color}
\usepackage{amssymb}
\usepackage{amssymb}
\usepackage[T1]{fontenc}
\usepackage{latexsym}
\usepackage{xypic}
\usepackage{eufrak}
\usepackage{euscript}
\usepackage{amsfonts,amsmath}
\usepackage{verbatim}
\usepackage{fancyhdr}
\usepackage{mathrsfs}
\usepackage{units}

\usepackage{hyperref}

\newtheorem{prop}{Proposition}[section]
\newtheorem{cor}[prop]{Corollary}

\newtheorem{thm}[prop]{Theorem}
\newtheorem{theorem}[prop]{Theorem}

\renewcommand{\geq}{\geqslant}
\def\leq{\leqslant}

\newcommand{\N}{\mathbb{N}}
\newcommand{\Z}{\mathbb{Z}}

\newcommand{\R}{\mathbb{R}}

\def\HH{\EuFrak H}

\def\e{\varepsilon}

\def\1{{\mathbf{1}}}

\def\sk{{\mathbb{D}}}

\def\1{{\mathbf{1}}}
\def\0.5{{\frac{1}{2}}}

\newcommand{\fin}
{ \vspace{-0.6cm}
\begin{flushright}
\mbox{$\Box$}
\end{flushright}
\noindent }

\newcommand{\qed}{\nopagebreak\hspace*{\fill}
{\vrule width6pt height6ptdepth0pt}\par}

\title{Fisher Information and the Fourth Moment Theorem}

\author{Ivan Nourdin\footnote{I. Nourdin is partially supported by the ANR grant `Malliavin, Stein and Stochastic Equations with Irregular Coefficients' [ANR- 10-BLAN-0121].} \,\,and David Nualart\thanks{D. Nualart is supported by the NSF grant  DMS1208625}}

\begin{document}
\maketitle
\begin{abstract}
Using a representation of the score function by means of the divergence operator we exhibit a sufficient condition, in terms of the negative moments of the norm of the Malliavin derivative, under which 
convergence in Fisher information to the standard Gaussian of sequences 
belonging to a given Wiener chaos is actually equivalent to convergence of only the fourth moment. Thus, our result may be considered as a further building block associated to the recent but already rich literature dedicated to the Fourth Moment Theorem of Nualart and Peccati \cite{nualartpeccati}.
To illustrate the power of our approach we prove a local limit theorem together with some rates of convergence for the normal convergence of a standardized version of the quadratic variation of the fractional Brownian motion.

\end{abstract}
 \textbf{Keywords:} Fisher information; total variation distance; relative entropy; Fourth Moment Theorem; Fractional Brownian motion; Malliavin calculus.
 
\section {Introduction}

Measuring the discrepancy between the law of a given real-valued random variable $F$ and that of its {\it Gaussian} counterpart $N$ is arguably an important and recurrent problem both in probability and statistics. For instance, one faces this situation when trying to prove a central limit type theorem, or when wanting to check the asymptotic normality of an estimator.
And quite often, the choice of a suitable probability metric reveals to be a crucial step. 

In the present paper, we are concerned with this question within the framework of the Malliavin calculus. 
More precisely, we will focus on the Wiener chaos of a given order and, as a way to measure the proximity betweens laws, we will work either with the $L^r$-distance between densities (especially for $r=1$ and $r=\infty$), or with the relative entropy $D(F\|N)$, or with 
the relative Fisher information $J(F)-1$. These three notions, that we will recall now, are strongly related to each other.

Let $F$ be a centered real-valued random variable with unit variance and density $p_F$. We suppose throughout that all needed  assumptions on $p_F$ (such as its strictly positivity, differentiability, etc.) are  always satisfied when required. 
Let also $N\sim N(0,1)$ be standard Gaussian, with density $p_N(x)=e^{-x^2/2}/\sqrt{2\pi}$, $x\in\R$.

The {\it $L^r$-distance} between densities of $F$ and $N$ is given by
\begin{eqnarray}
\|p_{F}-p_N\|_r &=& \left(\int_{\R}|p_F(x)-p_N(x)|^rdx\right)^\frac1r,\quad r\in[1,\infty);\label{lp}\\
\|p_F-p_N\|_\infty\!\! &=& {\rm sup}_{x\in\R}|p_F(x)-p_N(x)|\quad
\mbox{(assuming, say, that $p_F$ is continous)}.\notag
\end{eqnarray}
Actually, in what follows we will only consider the particular cases $r=1$ and $r=\infty$. This is because the bounds we will produce are going to be of the same order. So, a bound for the $L^r$-distance will simply follow from the crude estimate:
\[
\|p_{F}-p_N\|_r\leq \|p_F-p_N\|_1^{1/r}\,\, \|p_F-p_N\|_\infty^{1-1/r}.
\]

When $r=1$ in (\ref{lp}), it is an easy exercise (sometimes referred to as the Scheff\'e's theorem) to show that
$\|p_{F}-p_N\|_1 = 2 d_{TV}(F,N)$,
where $d_{TV}(F,N)$ is the {\it total variation distance} defined as
\begin{equation}\label{dtv}
d_{TV}(F,N)=\sup_{A\in\mathcal{B}(\R)}|P(F\in A)-P(N\in A)|.
\end{equation} 
It is clear from its very definition (\ref{dtv}) that
$d_{TV}(F,N)$ represents a strong measure on how close the laws of $F$ and $N$ are.

The {\it relative entropy} $D(F\|N)$ of $F$ with respect to $N$ is given by
\begin{equation}\label{relative}
D(F\|N)=\int_\R p_{F}(x)\log(p_{F}(x)/p_N(x))dx.
\end{equation}
Our interest in this quantity comes from its link with the total variation distance, as provided by the celebrated Csisz\'ar-Kullback-Pinsker inequality, according to which:
\begin{equation}\label{pinsker}
2\big(d_{TV}(F,N)\big)^2\leq D(F\|N).
\end{equation}
(In particular, note that $D(F\|N)\geq 0$.) See, e.g., \cite{bolleyvillani} for a proof of (\ref{pinsker}) and original references.

Inequality (\ref{pinsker}) shows that bounds on the relative entropy  translate directly into bounds on the total variation distance. Hence, it makes perfectly sense to quantify the
discrepancy between the law of $F$ and that of the standard Gaussian $N$ in terms of its relative entropy.
Actually, one can go even further by considering the {\it Fisher information} $J(F)$ of $F$. Let us recall its definition. Let $s_{F}(F)$ denote the {\it score} associated to $F$. This is the $F$-measurable random variable uniquely determined by
the following integration by parts:
\begin{equation}\label{scoredef}
E[\phi'(F)]=-E[s_F(F)\phi(F)]\quad \mbox{for all test function $\phi:\R\to\R$}.
\end{equation}
When it makes sense, it is easy to compute that $s_{F}=p_{F}'/p_F$.
Set $J(F)=E[s_{F}(F)^2]$ if the random variable $s_F(F)$ is square-integrable and $J(F)=+\infty$ otherwise. In the former case, it is a straightforward exercise to check that 
\[
J(F)-1=E[(s_F(F)+F)^2].
\]
In particular, $J(F)\geq 1=J(N)$ with equality if and only if $F$ is standard Gaussian. Our interest in the relative Fisher information $J(F)-1$ comes from its link with the relative entropy through the following de Bruijn's formula (stated in an integral and rescaled version due to Barron \cite{barron}; see also \cite[Theorem C.1]{Johnson}). Assume, without loss of generality, that $F$ and $N$ are independent; then
\begin{equation}\label{debruijn}
D(F\|N)=\int_0^1 \frac{J(\sqrt{t}F+\sqrt{1-t}N)-1}{2t}\,dt.
\end{equation}
Since from, e.g.,  \cite[Lemma 1.21]{Johnson} one has
$J(\sqrt{t}F+\sqrt{1-t}N)\leq tJ(F)+(1-t)J(N)=1+t(J(F)-1)$,
we deduce that
\begin{equation}\label{logsob}
D(F\|N)\leq \frac12 (J(F)-1).
\end{equation}
By comparing (\ref{logsob}) with (\ref{pinsker}), we observe that the gap between $J(F)$ and $1=J(N)$ is an even stronger measure of how the law of $F$ is close to the standard Gaussian $N$.
This claim is even more supported by the Shimizu's inequality \cite{shimizu}, which
gives a $L^\infty$-bound between $p_F$ and $p_N$ provided $p_F$ is continuous and satisfies $x^2p_F(x)\to 0$ as $x\to\pm\infty$:
\begin{equation}\label{shimizu}
\|p_F-p_N\|_\infty\leq \sqrt{J(F)-1}.
\end{equation}
(In the original statement of Shimizu \cite{shimizu}, there is actually an extra factor $(1+\sqrt{6/\pi})$ in the right-hand side of (\ref{shimizu}); but this latter was removed by Ley and Swan in \cite{leyswan}).\\

Let us now come to the description of the main results contained in the present paper.
From now on, we will systematically assume that $F$  
belongs to a Wiener chaos $\mathcal{H}_q$ of order $q\geq 2$,
that is, has the form of a $q$th multiple Wiener-It\^o integral
(see Section 2 below for precise definitions). 
Our first result is the following, with $\|DF\|$ the norm of the Malliavin derivative of $F$ (again, see Section 2 for details).
\begin{thm}\label{main1}
Let $q\geq 2$ be an integer and let $F\in\mathcal{H}_q$ have
variance one. Assume in addition that $\e>0$ and $\eta\geq 1$ satisfy
\begin{equation}\label{negative}
E[\|DF\|^{-4-\e}]\leq \eta.
\end{equation}
Then, there exists a constant $c>0$, depending on $q$, $\e$ and $\eta$ but \underline{not} on $F$, such that
\begin{equation}\label{fisher}
J(F)-1\leq c(E[F^4]-3).
\end{equation}
\end{thm}

In the next result, we take advantage of the conclusion (\ref{fisher}) of Theorem \ref{main1} to complete the current state of the art related to the {\it Fourth Moment Theorem} of Nualart and Peccati \cite{nualartpeccati}. See also the discussion located just after the statement of Corollary \ref{main2}.

\begin{cor}\label{main2}
Fix an integer $q\geq 2$, and let $(F_n)\subset \mathcal{H}_q$  be a sequence of random variables satisfying $E[F_n^2]=1$ for all $n$.
Then, the following four assertions are equivalent
as $n\to\infty$:
\begin{enumerate}
\item[(a)] $E[F_n^4]\to 3$;
\item[(b)] $F_n\overset{\rm law}{\to} N\sim N(0,1)$;
\item[(c)] $d_{TV}(F_n,N)=\frac12\|p_{F_n}-p_N\|_1\to 0$;
\item[(d)] $D(F_n\|N)\to 0$.
\end{enumerate}
Moreover, there exists $c_1,c_2,c_3,c_4>0$ (independent of $n$) such that,
for all $n$ large enough,
\begin{eqnarray}
c_1\max\{|E[F_n^3]|,E[F_n^4]-3\}\leq d_{TV}(F_n,N)&\leq& c_2\max\{|E[F_n^3]|,E[F_n^4]-3\}\label{ineq1}\\
&\leq&c_3\sqrt{E[F_n^4]-3}\label{ineq1bis};\\
\notag\\
D(F_n\|N)&\leq&c_4 (E[F_n^4]-3)|\log(E[F_n^4]-3)|.\label{ineq2}
\end{eqnarray}
Suppose in addition that, for some $\e>0$,
\begin{equation}\label{H}
\limsup_{n\to\infty} E[\|DF_n\|^{-4-\e}]<\infty.
\end{equation}
Then, the four previous assertions $(a)-(d)$ are equivalent to the following two further assertions:
\begin{enumerate}
\item[(e)] $\|p_{F_n}-p_N\|_\infty\to 0$;
\item[(f)] $J(F_n)\to 1$.
\end{enumerate}
More precisely, one has the existence of $c_5>0$ (independent of $n$) such that, for all $n$ large enough,
\begin{eqnarray}
\|p_{F_n}-p_N\|_\infty^2&\leq& c_5(E[F_n^4]-3)\label{david}\\
J(F_n)-1
&\leq& c_5(E[F_n^4]-3).\label{boundcor}
\end{eqnarray}
\end{cor}

\bigskip

Equivalence between $(a)$ and $(b)$ in Corollary \ref{main2} is known as the {\it Fourth Moment Theorem}. This striking result, discovered by Nualart and Peccati in \cite{nualartpeccati}, has been the starting point of a new and fruitful line of research, consisting in using the Malliavin calculus to prove limit theorems. It has led to a burst of new research in many different fields, such as information theory \cite{NPS}, stochastic geometry \cite{lachiezepeccati,reitznerschulte}, Markov operator \cite{ACP,ledoux}, random matrices of large size \cite{NPalea}, free probability \cite{KNPS,NPpoisson}, $q$-calculus \cite{arizmendi,DNN}, computer science \cite{de1,de2}, cosmology \cite{marinuccipeccati,marinucciwigman}, statistics \cite{bardetsurgailis,hunualart}, or spin glasses \cite{NPV,viens}, to name a few.
One can also consult the constantly updated webpage\\

\url{http://www.iecn.u-nancy.fr/~nourdin/steinmalliavin.htm}\\
\\
for literally hundreds of results related to the Fourth Moment Theorem and its ramifications.

Equivalence between $(a)$ and $(c)$  in Corollary \ref{main2}, coming from the bound 
(\ref{ineq1bis}), is due to Nourdin and Peccati \cite{NourdinPeccatiPTRF}.
By combining Malliavin calculus with the Stein's method, they were indeed able to show that, 
for any random variable $F\in\mathcal{H}_q$ such that $E[F^2]=1$,
\begin{equation}\label{NP}
d_{TV}(F,N)\leq 
2E|1-\frac1q\|DF\|^2|\leq
\sqrt{\frac{4q-4}{3q}}\sqrt{E[F^4]-3},
\end{equation}
see \cite[Theorem 5.2.6]{nourdinpeccatibook}. 
(Note that $E[F^4]>3$ in (\ref{NP}): see, e.g., \cite[Lemma 5.2.4]{nourdinpeccatibook}.)
The refinement (\ref{ineq1}) of (\ref{ineq1bis}), leading to optimal rates, is taken from \cite{best}.

Equivalence between $(a)$ and $(d)$  in Corollary \ref{main2}, as well as the
bound (\ref{ineq2}), was shown by Nourdin, Peccati and Swan in \cite{NPS}. 
Their strategy of proof relied on the discovery a novel representation formula for the relative entropy, namely,
\[
D(F\|N)=\frac12\int_0^1 \frac{t}{1-t}E\big[E[N(1-\frac1q\|DF\|^2)|\sqrt{t}F+\sqrt{1-t}N]^2\big]dt
\]
for any $F\in\mathcal{H}_q$ with unit variance and where $N\sim N(0,1)$
is supposed to be independent of $F$.

When (\ref{H}) is satisfied for $\e=2$, the inequality 
(\ref{david})
(leading to the equivalence between $(a)$ and $(e)$ in Corollary \ref{main2})
was proved by Hu, Lu and Nualart in \cite{HLN}, after adapting Stein's method to handle the supremum distance.
Note that combining our Theorem \ref{main1} with Shimizu inequality (\ref{shimizu}) allows to recover (\ref{david}) (which corresponds to Theorem 4.1 in \cite{HLN}).

Finally, inequality (\ref{boundcor}) (leading to the equivalence between $(a)$ and $(f)$ in Corollary \ref{main2})  is new and will be a direct consequence of Theorem \ref{main1}.
It is worth noting at this stage that validity of (\ref{negative}) is, unfortunately, far to be a small assumption.
Let us discuss this point a little bit more.
As it is well-known, the Bouleau-Hirsch criterion (see, e.g., \cite[Theorem 2.1.3]{nualartbook}) asserts that any (smooth and bounded enough) random variable $F$ in the Wiener space admits a density as soon as $P(\|DF\|>0)=1$. This latter condition is always satisfied
for $F\in\mathcal{H}_q$ with unit variance, see  \cite{Shigekawa}. In fact, one can prove a far better statement, see \cite[(3.19)]{nourdinpoly}: there exists a constant $c_q>0$ such that, for all $x>0$ and all $F\in\mathcal{H}_q$ with unit variance,
\begin{eqnarray}\label{CR}
P(\|DF\|^2\geq x)&\ge& 1-c_q\,x^{\frac{1}{2q-2}}.
\end{eqnarray}
As a consequence, using that 
\begin{equation}
E[\|DF\|^{-r}]=\int_0^\infty P(\|DF\|^{2}\leq u^{-\frac{2}{r}})du
\leq 1+\int_1^\infty P(\|DF\|^{2}\leq u^{-\frac{2}{r}})du,\label{118bis}
\end{equation}
one deduces from (\ref{CR}) that, for all $F\in\mathcal{H}_q$ with unit variance,
\begin{equation}\label{car}
E[\|DF\|^{-r}]\leq 1+\frac{c_q}{\frac1{r(q-1)}-1} \quad\mbox{provided $r<\frac{1}{q-1}$}.
\end{equation}

Unfortunately, one cannot deduce ({\ref{H}) from (\ref{car}). It means that verifying (\ref{H}) has to be made on a case-by-case basis, and heavily depends  on the particular sequence $(F_n)$ we are dealing with. In \cite{HLN}, one can find an application for the least squares estimator of the parameter $\theta$ in the Ornstein-Uhlenbeck process $dX_t = −\theta X_tdt + dW_t$, where $W$ is a standard Brownian motion.
In the present paper, we consider a more involved application to the quadratic variation of a fractional Brownian motion $B^H$ of index $H$. We obtain optimal rates for the relative Fisher information when $H<\frac58$, and (possibly suboptimal) rates when $H<\frac34$.
More precisely,
let us introduce the so-called {\it fractional Gaussian noise} associated with $B^H$, which is the Gaussian sequence given by
\begin{equation}\label{xk}
\xi_k=B^{H}(k+1)-B^H(k),\quad k\in\N\cup\{0\}.
\end{equation}
Set
\begin{equation}\label{fn}
F_n := \frac{1}{\sqrt{n\,v_n}}\sum_{k=0}^{n-1} (\xi_k^2-1),
\end{equation}
with $v_n>0$ chosen so that $E[F_n^2]=1$. 
It is well-known (it is indeed a very particular case of the Breuer-Major Theorem \cite{breuermajor}, see also \cite[Theorem 7.2]{nourdinbookfbm}) that, as $n\to\infty$,
\begin{equation}\label{cvlaw}
F_n\overset{\rm law}{\to}N(0,1)\quad\Longleftrightarrow \quad H\in(0,3/4].
\end{equation}
In Section 4, we will show that (\ref{H}) is satisfied for $F_n$ defined by (\ref{fn}). Then, as a consequence of (\ref{david}) and (\ref{boundcor}) on one hand and of the estimates
for $E[F_n^4]-3$ computed in \cite{BBNP} on the other hand, we will be able to deduce
the following {\it local} limit theorem for $F_n$.
\begin{theorem}\label{fractional}
Let $F_n$ be as in (\ref{fn}). Then, 
their exists $c,C>0$ independent of $n$ such that, for all $n$ large enough,
\[
\|p_{F_n}-p_N\|_\infty \leq \sqrt{J(F_n)-1}\leq C\,\times
\left\{
\begin{array}{lll}
n^{-\frac12}&\mbox{ if }&0<H<\frac58\\
n^{-\frac12}\log^{\frac32}n&\mbox{ if }&H=\frac58\\
n^{4H-3}&\mbox{ if }&\frac58<H<\frac34
\end{array}
\right.
\]
and
\[
\sqrt{J(F_n)-1}\geq c\,n^{-\frac12}\quad\mbox{if $H<\frac58$}.
\]
\end{theorem}

\bigskip

 A brief outline of the paper is as follows.
   In Section 2, we introduce the language of the Malliavin calculus, which is the framework in which our study takes place. We also recall the Carbery-Wright inequality, which will play a key role in the proof of Theorem \ref{fractional}. Proofs of Theorem \ref{main1} and Corollary \ref{main2} are presented in Section 3. Finally, Section 4 contains the proof of Theorem \ref{fractional}.

\section{Notation and preliminaries}

\subsection{The language of Gaussian analysis and Malliavin calculus}\label{ss:mall}

We start by briefly recalling some basic notation and results connected to Gaussian analysis and Malliavin calculus. The reader is referred to \cite{nourdinpeccatibook,nualartbook} for details or missing proofs.

\medskip

Let $\HH$ be a real separable Hilbert space 
with inner product  $\langle\cdot,\cdot\rangle_\HH$.  
The norm of $\HH$ will be denoted  by $\|\cdot \|=\| \cdot \|_\HH$.
 Recall 
that we call
{\it isonormal Gaussian process} 
over 
$\HH$ 
any
centered 
Gaussian 
family $X = \{X(h) : h\in \HH\}$, defined on a probability space $(\Omega, \mathcal{F}, P)$ and such that $E[X(h) X(g)] = \langle h, g\rangle_\HH$ for every $h,g\in \HH$. Assume from now on that $\mathcal{F}$ is
the $\sigma$-field
generated by $X$.

\medskip

For any integer $q\in\N\cup\{0\}$, we denote by $\mathcal{H}_q$ the $q$th {\it Wiener chaos} of $X$. 
We recall that $\mathcal{H}_0$ is simply $\R$ whereas, for any $q\geq 1$, $\mathcal{H}_q$ is the closed linear subspace of $L^2(\Omega)$ generated by the family of random variables $\{H_q(X(h)), h\in \mathfrak{H}, \|h\|_{\mathfrak{H}}=1\}$, with $H_q$ the $q$th Hermite polynomial given by
\[
H_q(x)= (-1)^q e^{\frac {x^2}2} \frac {d^q} {dx^q} \left( e^{-\frac {x^2}2}\right).
\]
For any $q\geq 1$, 
we denote by $\mathfrak{H}^{\otimes q}$ (resp. $\mathfrak{H}^{\odot q}$) the $q$th tensor product (resp. the $q$th {\it symmetric} tensor product) of $\mathfrak{H}$. Then, the mapping 
$I_q(h^{\otimes q})= H_q(X(h))$ can be extended to a linear isometry between  $\mathfrak{H}^{\odot q}$ (equipped with the modified norm $\sqrt{q!} \| \cdot \|_{\mathfrak{H}^{\otimes q}}$) and $\mathcal{H}_q$. 
For $q=0$ and $x\in \mathbb{R}$, we write $I_0(x)=x$.
In the
particular case where $\HH=L^2(A,\mathscr{A},\mu)$, where
$(A,\mathscr{A})$ is a measurable space and $\mu$ is a
$\sigma$-finite and non-atomic measure, one has that $\HH^{\odot
q}= L^2_s(A^q,\mathscr{A}^{\otimes q}, \mu^{\otimes q})$ is
the space of symmetric and square integrable functions on $A^q$.
Moreover, for every $f\in\HH^{\odot q}$, the random variable
$I_q(f)$ coincides with the multiple Wiener-It\^o integral (of
order $q$) of $f$ with respect to $X$.
\medskip

Recall that $L^2(\Omega)= \bigoplus_{q=0}^\infty \mathcal{H}_q$, meaning that every square-integrable random variable $F$ measurable with respect to $\mathcal{F}$ admits a unique decomposition of the type
\begin{equation}\label{e:hh21}
F = E[F] +\sum_{q=1}^\infty I_q(f_q),
\end{equation}
where the series converges in $L^2(\Omega)$, and $f_q \in \HH^{\odot q}$, for $q\geq 1$. Identity (\ref{e:hh21}) is the so-called {\it Wiener-It\^o chaotic decomposition} of $F$.
According to a classical result of Shigekawa \cite{Shigekawa}, when $F$ is not zero and when the kernels $f_q$ in \eqref{e:hh21} all equal zero except for a finite number, then the distribution of $F$ necessarily admits a density with respect to the Lebesgue measure.

\medskip

Let $\{e_i, i\geq 1\}$ be a complete orthonormal system in $\mathfrak{H}$. Given $f\in \mathfrak{H}^{\odot p}$ and $g\in \mathfrak{H}^{\odot q}$, for every $r=0,\dots, p\wedge q$, the \textit{contraction} of $f$ and $g$ of order $r$ is the element of $\mathfrak{H}^{\otimes(p+q-2r)}$ defined by
\[
f\otimes_r g= \sum_{i_1, \dots, i_r=1}^\infty \langle f, e_{i_1} \otimes \cdots \otimes e_{i_r}
\rangle_{\mathfrak{H}^{\otimes r}} \otimes  \langle g, e_{i_1} \otimes \cdots \otimes e_{i_r}
\rangle_{\mathfrak{H}^{\otimes r}}.
\]
Note that, in the particular case where $\HH=L^2(A,\mathscr{A},\mu)$ (with $\mu$ non-atomic), one has that
\begin{eqnarray*}
&&(f\otimes_r g)(t_1,\ldots,t_{p+q-2r}) \\
&=& \int_{A^r}
f(t_1,\ldots,t_{p-r},s_1,\ldots,s_r)\,g(t_{p-r+1},\ldots,
t_{p+q-2r},s_1,\ldots,s_r)d\mu(s_1)\ldots d\mu(s_r).
\end{eqnarray*}
Moreover, $f\otimes_0 g=f\otimes g$ equals the tensor product of
$f$ and $g$ while, for $p=q$, $f\otimes_p g=\langle
f,g\rangle_{\HH^{\otimes p}}$.
The contraction $f\otimes_r g$ is not necessarily symmetric, and we denote by $f \widetilde{\otimes}_rg$ its symmetrization.
We have the following product formula: if $f\in \EuFrak
H^{\odot p}$ and $g\in \EuFrak H^{\odot q}$ then
\begin{eqnarray}\label{multiplication}
I_p(f) I_q(g) = \sum_{r=0}^{p \wedge q} r! {p \choose r}{q \choose r} I_{p+q-2r} (f\widetilde{\otimes}_{r}g).
\end{eqnarray}

\medskip

We will now introduce some standard operators from Malliavin calculus. 
Let $\mathcal{S}$ be the set of all cylindrical random variables of the form
\[
F=g(X(h_1), \dots, X(h_n)),
\]
where $n\geq 1$, $h_i \in \mathfrak{H}$, and $g$ is infinitely differentiable such that all its partial derivatives have polynomial growth. The Malliavin derivative of $F$ is the element of $L^2(\Omega;\mathfrak{H})$ defined by
\[
DF= \sum_{i=1}^n \frac {\partial g}{\partial x_i}(X(h_1), \dots, X(h_n)) h_i.
\]
By iteration, for every $m\geq 2$, we define the $m$th derivative $D^mF$ which is an element of $L^2(\Omega; \mathfrak{H}^{\odot m})$. For $m\geq 1$ and $p\geq 1$, $\mathbb{D}^{m,p}$ denote the closure of $\mathcal{S}$ with respect to the norm $\| \cdot \|_{m,p}$ defined by
\[
\|F\|^p_{m,p} = E[|F|^p] + \sum_{j=1}^m E\left[ \|D^jF\|^p_{\mathfrak{H}^{\otimes j}}\right].
\]
One can then extend the definition of $D^m$ to $\mathbb{D}^{m,p}$. 
When $m=1$, one simply write $D$ instead of $D^1$. 
As a consequence of the hypercontractivity property of the Ornstein-Uhlenbeck semigroup (see, e.g., \cite[Theorem 2.7.2]{nourdinpeccatibook}),
all the $\| \cdot \|_{m,p}$-norms are equivalent in any {\it finite} sum of Wiener chaoses. This is a crucial result that will be used all along the paper.

The
Malliavin derivative $D$ satisfies the following \textsl{chain
rule}: if $\varphi:\R^n\rightarrow\R$ is in $\mathcal{C}^1_b$
(that is, belongs to the set of continuously differentiable
functions with a bounded derivative) and if
$\{F_i\}_{i=1,\ldots,n}$ is a vector of elements of $\sk^{1,2}$,
then $\varphi(F_1,\ldots,F_n)\in\sk^{1,2}$ and
$$
D\varphi(F_1,\ldots,F_n)=\sum_{i=1}^n
\frac{\partial\varphi}{\partial x_i} (F_1,\ldots, F_n)DF_i.
$$
Also, when $\HH=L^2(A,\mathscr{A},\mu)$ (with $\mu$ non-atomic), one has, for any $f\in L^2_s(A^q,\mathscr{A}^{\otimes q}, \mu^{\otimes q})$,
\[
D_x(I_q(f))=qI_{q-1}(f(\cdot,x)),\quad x\in A.
\]

\medskip 

The {\it divergence operator} $\delta$, which will play a crucial role in our approach,  
is defined as the adjoint of $D$. 
Denoting by ${\rm dom}\, \delta$ 
its domain, 
one 
has 
the 
so-called 
{\it integration by parts formula}: 
for 
every 
$D\in \mathbb{D}^{1,2}$ 
and 
every $u\in {\rm dom}\delta$,
\begin{equation}\label{e:ipp21}
E[F\delta(u)] = E[\langle DF ,u\rangle_\HH].
\end{equation}
We will moreover need the following two properties. For every $F\in\mathbb{D}^{1,2}$ and every $u\in {\rm dom}\delta$ such
that $Fu$ and $F\delta(u)+\langle DF,u\rangle_\HH$ are square integrable, one
has that $Fu\in{\rm dom}\delta$ and
\begin{equation}\label{factorout}
\delta(Fu)=F\delta(u)-\langle DF,u\rangle_\HH.
\end{equation}
Also, one has a commutation relationship between the
Malliavin derivative and the Skorohod integral:
\begin{equation}
D\delta (u)=u+\delta (Du),  \label{comm1}
\end{equation}%
for any $u\in \mathbb{D}^{2,2}(\EuFrak H)$. 
In particular, for such an $u$,
\begin{equation}\label{moment2}
E[\delta(u)^2]=E[\|u\|_\HH^2]+E[\|Du\|^2_{\HH^{\otimes 2}}].
\end{equation}

\subsection{Carbery-Wright inequality}

In the proof of Theorem \ref{fractional}, we will make use of the following inequality due to Carbery and Wright \cite[Theorem 8]{CW}: there is an absolute constant $c>0$ such that, for any $d,n\geq 1$, any polynomial $Q:\mathbb{R}^n \rightarrow \mathbb{R}$ of degree at most $d$ and  any Gaussian random vector $(X_1,\ldots,X_n)$,
\begin{equation}  \label{ff1}
E[|Q(X_1, \dots, X_n)|] ^{\frac 1{d}} P(|Q(X_1, \dots, X_n)| \leq x) \leq c\,d\,x^{\frac 1d},\quad x>0.
\end{equation}

\section{Proofs of Theorem \ref{main1} and Corollary \ref{main2}}

In what follows, 
$c$ denote positive constants which may depend of $q$, $\e$ and $\eta$ but \underline{not} of $F$, and whose values may change from one appearance to the next. Also, $\langle\cdot,\cdot\rangle$ ($\|\cdot\|$, respectively) always stands for inner product (the norm, respectively) in an appropriate tensor product $\HH^{\otimes s}$.

\subsection{Proof of Theorem \ref{main1}}

Observe
first that, without loss of generality, we may and will assume that $X$ is
an isonormal process over some Hilbert space of the type $\HH=L^2(A,\mathscr{A},\mu)$ (with $\mu$ non atomic).

Due to (\ref{negative}) and the fact that $F$ has moments of all order by hypercontractivity, it is straightforward to check that $DF\|DF\|^{-2}\in {\rm dom}\delta$ with
$E[\delta(DF\|DF\|^{-2})^2]<\infty$. 
Let $\phi:\R\to\R$ be a test function. 
We have, on one hand,
\[
E[\delta(DF\|DF\|^{-2})\phi(F)]=E[\langle D\phi(F),DF\|DF\|^{-2}\rangle]
=E[\phi'(F)].
\]
After setting  $\Sigma = 1-\frac1q \|DF\|^2$ and because $\delta DF= qF$ for any $F\in\mathcal{H}_q$, we deduce from (\ref{scoredef}) that
\[
s_F(F)+F=-E\left[\delta\left(DF\big(\|DF\|^{-2}-\frac1q\big)\right)\bigg|F\right]
=-E\big[\delta(DF\|DF\|^{-2}\Sigma)\big| F\big].
\]
Using the formula 
 \[
 \delta(GDF)=  G \delta(DF)- \langle DF,DG\rangle=q FG -\langle DF,DG\rangle,
 \]
 one can write for  $G= \|DF\|^{-2}  \Sigma$,
 \[
 -\delta ( DF \|DF\|^{-2}  \Sigma)= -q F\|DF\|^{-2}  \Sigma + \langle DF, DG\rangle.
 \]
 Notice that
 \[
 DG= -2 \|DF\|^{-4}  (D^2F \otimes_1 DF)\Sigma -\frac 2q   \|DF\|^{-2}  (D^2F \otimes_1 DF).
 \]
 Therefore,
 \[
 \langle DF, DG\rangle=  -2 \|DF\|^{-4}  \langle D^2F ,  DF\otimes DF )\Sigma -\frac 2q   \|DF\|^{-2}  \langle DF, D^2F \otimes_1 DF\rangle.
 \]
 This leads to the estimate
 \[
 |\langle DF, DG\rangle| \le  2 \|DF\|^{-2}  \|D^2F\| |\Sigma| + \frac 2q  \|DF\|^{-1}  \|D^2F \otimes_1 DF\|.
 \]
 As a consequence,
 \[
 \delta ( DF \|DF\|^{-2}  \Sigma)^2 \le 2 q^2 F^2\|DF\|^{-4}  \Sigma^2+ 16 \|DF\|^{-4}  \|D^2F\|^2 \Sigma^2
 + \frac {16}{q^2} \|DF\|^{-2}  \|D^2F \otimes_1 DF\|^2.
 \]
Thus, using among other properties the hypercontractivity for $F$, $\|D^2F\|^2$ and $\Sigma$,
\begin{eqnarray}
J(F)-1
&=&E[(s_F(F)+F)^2]\leq
E\left[
 \delta ( DF \|DF\|^{-2}  \Sigma)^2 
\right]\notag\\
&\leq&
c\, E[\|DF\|^{-4-\e}]^{\frac4{4+\e}}\big(E[\Sigma^2]+E[\Sigma^2]E[\|D^2F\|^4]^\frac12+E[\|D^2F \otimes_1 DF\|^4]^\frac12\big).\notag
\\\label{eq07}
\end{eqnarray}
Now, use the product formula to get that, for any $x\in A$,
\begin{eqnarray*}
(D^2F\otimes_1 DF)(x)&=&q^2(q-1)\int_A I_{q-2}(f(x,y,\cdot))I_{q-1}(f(y,\cdot))d\mu(y)\\
&=&q^2(q-1)\sum_{r=0}^{q-2}r!\binom{q-1}{r}\binom{q-2}{r}I_{2q-3-2r}\left(\int_A f(x,y,\cdot)\otimes_r f(y,\cdot)d\mu(y)\right)\\
&=&q^2(q-1)\sum_{r=1}^{q-1} (r-1)!\binom{q-1}{r-1}
\binom{q-2}{r-1}I_{2q-1-2r}((f\otimes_r f)(x,\cdot)).
\end{eqnarray*}
As a result, using again the product formula and with $c_{q,r,s,a}$ some constant whose exact value is useless here,
\begin{eqnarray*}
&&\|D^2F\otimes_1 DF\|^2\\
&=&\sum_{r,s=1}^{q-1}\sum_{a=0}^{2q-1-2(r\vee s)}
c_{q,r,s,a}\,I_{4q-2-2r-2s-2a}\left(
\int_A \left(
\widetilde{(f\otimes_r f)(x,\cdot)}\otimes_a\widetilde{(f\otimes_s f)(x,\cdot)}
\right)d\mu(x)
\right),
\end{eqnarray*}
implying in turn
\begin{eqnarray*}
&&E\|D^2F\otimes_1 DF\|^4\\
&\leq& c\sum_{r,s=1}^{q-1} \sum_{a=0}^{2q-1-2(r\vee s)}
\left\| \int_A \left(
\widetilde{(f\otimes_r f)(x,\cdot)}\widetilde{\otimes}_a\widetilde{(f\otimes_s f)(x,\cdot)}
\right)d\mu(x)\right\|^2\\
&\leq& c\sum_{r,s=1}^{q-1} \sum_{a=0}^{2q-1-2(r\vee s)}
\left( \int_A \left\|
\widetilde{(f\otimes_r f)(x,\cdot)}\widetilde{\otimes}_a\widetilde{(f\otimes_s f)(x,\cdot)}
\right\|d\mu(x)\right)^2\\
&\leq&c\sum_{r,s=1}^{q-1}
\left(
\int_A \|(f\otimes_r f)(x,\cdot)\|\|(f\otimes_s f)(x,\cdot)\|d\mu(x)
\right)^2\\
&\leq& c \left(
\sum_{r=1}^{q-1}\int_A \|(f\otimes_r f)(x,\cdot)\|^2d\mu(x)
\right)^2=c\left(
\sum_{r=1}^{q-1}\|f\otimes_r f\|^2
\right)^2\leq c(E[F^4]-3)^2,
\end{eqnarray*}
the last inequality following from \cite[identities (5.2.5)-(5.2.6)]{nourdinpeccatibook}. On the other hand, we  can also write
\begin{eqnarray*}
\|D^2F\|^2&=&q^2(q-1)^2 \int_{A^2} I_{q-2}(f(x,y,\cdot))^2d\mu(x)d\mu(y)\\
&=&q^2(q-1)^2\sum_{r=0}^{q-2}r!\binom{q-2}{r}^2I_{2q-4-2r}
\left(\int_{A^2}f(x,y,\cdot)\otimes_r f(x,y,\cdot)d\mu(x)d\mu(y)\right)\\
&=&q^2(q-1)^2\sum_{r=2}^q (r-2)!\binom{q-2}{r-2}^2 I_{2q-2r}(f\otimes_rf),
\end{eqnarray*}
so that, using moreover that $\|f\widetilde{\otimes}_r f\|^2\leq \|f\|^4 =1/q!^2$,
\[
E[\|D^2F\|^4]=q^4(q-1)^4 \sum_{r=2}^q (r-2)!^2\binom{q-2}{r-2}^4(2q-2r)!\|f\widetilde{\otimes}_r f \|^2\leq c.
\]
Finally, recall from \cite[Lemma 5.2.4]{nourdinpeccatibook} that $E[\Sigma^2]\leq\frac{q-1}{3q}(E[F^4]-3)$. Hence,
by plugging all the previous estimates in (\ref{eq07}), one finally obtains the desired inequality (\ref{fisher}).
\fin

\subsection{Proof of Corollary \ref{main2}}

As we said in the Introduction, the equivalences between $(a)$, $(b)$, $(c)$, $(d)$ and $(e)$ (provided (\ref{negative}) holds true for the latter one), as well as the estimates (\ref{ineq1}) and (\ref{ineq2}), are straightforward consequences of the main results contained in \cite{HLN}, \cite{NourdinPeccatiPTRF}, \cite{NPS} and \cite{nualartpeccati}.

Note in passing that Shimizu's inequality (\ref{shimizu}) indeed takes place for the random variable we are considering. This is because, if $F\in\mathcal{H}_q$ satisfies $E[\|DF\|^{-4-\e}]<\infty$ then, by \cite[Proposition 2.1.1]{nualartbook}, $F$ has a continuous density given by
\[
p_F(x)=E[{\bf 1}_{\{F>x\}}\delta(DF\|DF\|^{-2})].
\]
Hence, using moreover that $E[\delta(DF\|DF\|^{-2})]=0$, one deduces that
\[
p_F^2(x)\leq P(|F|>x)E[\delta(DF\|DF\|^{-2})^2]\leq |x|^{-r}E\big[|F|^r\big]
E[\delta(DF\|DF\|^{-2})^2],
\]
implying in turn that $p(x)=o(|x|^{-n})$ as $|x|\to\infty$.

The equivalence}between $(f)$ and $(a)$, provided condition (\ref{H}) is fulfilled, follows immediately from Theorem \ref{main1} (for one implication) and the fact that bounds on the Fisher information translates directly into bounds on the total variation distance through (\ref{fisher}) and (\ref{pinsker}) (for the other implication).

Inequality (\ref{boundcor}) also follows immediately from (\ref{fisher}).\fin

\section{Proof of Theorem \ref{fractional}}

\subsection{Preparation to the proof}
For $F_n$ as in the statement, recall from \cite{BBNP,best} that
their exist $c,C>0$ independent of $n$ such that, for all $n$ large enough,
\[
E[F_n^4]-3\leq C\,\times
\left\{
\begin{array}{lll}
n^{-1}&\mbox{ if }&0<H<\frac58\\
n^{-1}\log^{3}n&\mbox{ if }&H=\frac58\\
n^{8H-6}&\mbox{ if }&\frac58<H<\frac34
\end{array}
\right.
\]
and
\[
d_{TV}(F_n,N)\geq c\,n^{-1}\quad \mbox{if $H>\frac58$}.
\]
Assume for an instant that (\ref{H}) has been checked. Then, as far as the upper (resp. lower) bound is concerned,
the desired conclusion directly follows from (\ref{david}) and (\ref{boundcor}) (from (\ref{pinsker}) and (\ref{logsob}), respectively).

So, in order to complete the proof of Theorem \ref{fractional}, it remains to check that (\ref{negative}) holds true.

\subsection{Checking (\ref{negative})}

For simplicity, throughout all the proof, we write $B$ instead of $B^H$ to indicate the fractional Brownian motion of index $H\in(0,1)$ we are dealing with. 
We know that $B$ has an integral representation of the form
\begin{equation}\label{repreK1}
B _t = \int_0^t K_1 (t,s)dW_s,\quad t\geq 0,
\end{equation}
where $\{W_t, t\ge 0\}$ is a standard Brownian motion. By convention we will assume that  $K_1 (t,s)=0$ if $t\le s$.  Also, we set $\Delta K_1(t,s)= K_1(t,s) -K_1(t-1,s)$.

Recall the definition (\ref{xk}) of $\xi_k$ and the definition  (\ref{fn}) of $F_n$.
We claim that, for any $p\ge 1$, there exist $n_0$ such that
\begin{equation}\label{claim}
\sup_{n\ge n_0} E(\|DF_n\|^{-p})<\infty.
\end{equation}
Note that (\ref{claim}) implies (\ref{H}). 
 
The proof of our claim (\ref{claim}) is based on the following approach. First we will derive a lower bound for  $\|DF_n\|^2$ (in distribution) denoted by  $B_n$ and defined in (\ref{Bn}), which can be expressed as the sum of the square norms of $n$  Gaussian random variables.  Then,  we fix an integer $N\ge 1$ and we decompose $B_n$ into the sum of $N$ blocks $B_n^i$, $i=0,\dots, N-1$  of size $[n/N]$. The basic inequality (\ref{ineq}) reduces the problem to estimate negative moments of order $\frac pN$ of each block $B^i_n$, provided these blocks are independent. These negative moments can be estimated by the Carbery-Wright inequality if $N$ is large enough. Actually, the blocks are not independent, but we can control the conditional expectation of each block given the previous ones, using the properties of the fractional Brownian motion. Then, it suffices to show that these conditional expectations  do not vanish as $n$ tends to infinity (condition (\ref{eq4})), which is done in the two final steps of the proof.

The random variables $\xi_k$ form a  centered stationary Gaussian sequence with covariance 
\[
\rho(k)= E[\xi_{r}\xi_{r+k}]= \frac 12 \left( |k+1|^{2H} + |k-1| ^{2H} - 2|k|^{2H} \right).
\]
We can thus write, with $D$ the Malliavin derivative with respect to $B$,
\begin{eqnarray*}
\|DF_n\|^2&=&  \frac 4{nv_n}\sum_{j,k=1}^n \xi_j \xi_k \rho(j-k).
\end{eqnarray*}
Suppose that $\{\widetilde{\xi}_j, 1\le j\le n\}$ is an independent copy of the sequence $\{\xi_j, 1\le j\le n\}$.
Then,
\[
\|DF_n\|^2 =\frac 4{nv_n}  \widetilde{E} \left( \left|  \sum_{j=1}^n \xi_j  \widetilde{\xi}_j \right|^2\right),
\]
with obvious notations.
The sequence $\widetilde{\xi}_j$ can be chosen of the form $\widetilde{\xi}_j= \widetilde{B}_j-\widetilde{B}_{j-1}$,
where  $\widetilde{B}$ is a fractional Brownian motion of Hurst index $H$ which is independent of $B$. We know that 
$\widetilde{B}$ has a representation of the form (different in nature from (\ref{repreK1}))
\[
\widetilde{B}_t= \kappa_H \left( \int_0^t (t-s) ^{H-\frac 12} d\widetilde{W}_s +Z_t \right),
\]
where $\widetilde{W}$ is a standard Brownian motion (independent of $W$), $Z$ is a process independent of $\widetilde{W}$ (and of $W$) and $\kappa_H$ is a constant only depending on $H$. We set $K_2(t)=t ^{H-\frac 12}$ if $t>0$ and $K_2(t)=0$ if $t\le 0$, and $\Delta K_2(t)= K_2(t) -K_2(t-1)$.
With this notation we can write
\[
\widetilde{\xi}_j= \kappa_H \left( \int_{0}^j \Delta K_2(j-s) d\widetilde{W}_s+ Z_j -Z_{j-1} \right).
\]
As a consequence, and since $v_n=\frac2n\sum_{k,l=1}^{n}\rho(k-l)^2\leq2\sum_{j\in\Z}\rho(j)^2<\infty$ for $H\in (0,\frac{3}{4})$,
\begin{eqnarray*}
\|DF_n\|^2& =&
\frac {c_H}n     \widetilde{E} \left( \left|  \sum_{j=1}^n \xi_j\int_{0}^j \Delta K_2(j-s) d\widetilde{W}_s \right|^2\right) \\
&=& \frac {c_H}n    \widetilde{E} \left( \left|  \sum_{j=1}^n \xi_j \sum_{h=1}^{j} \int_{h-1}^h  \Delta K_2(j-s)d\widetilde{W}_s
 \right|^2\right)\\
 &=& \frac {c_H}n   \widetilde{E} \left( \left|  \sum_{h=1}^n\int_{h-1}^h \left( \sum_{j=h}^{n}     \xi_j \Delta K_2(j-s)\right) d\widetilde{W}_s
 \right|^2\right) \\
 &=&  \frac {c_H}n      \sum_{h=1}^n\int_{h-1}^h \left( \sum_{j=h}^{n}     \xi_j \Delta K_2(j-s)\right) ^2ds.
\end{eqnarray*}
Making a change of indices, we obtain  
\[
\|DF_n\|^2 \ge   A_n :=\frac {c_H}n    \sum_{h=1}^n\int_{n-h}^{n-h+1} \left( \sum_{j=1}^{h}     \xi_{n-j+1} \Delta K_2(n-j+1-s)\right) ^2ds.
\]
The sequences $\{\xi_j, 1\le j\le n\}$ and $\{\xi_{n-j+1}, 1\le j\le n\}$ have the same law, so  $A_n $ has the same law as
\begin{equation} \label{Bn}
B_n:= \frac {c_H}n\sum_{h=1}^n\int_{n-h}^{n-h+1} \left( \sum_{j=1}^{h}     \xi_j \Delta K_2(n-j+1-s)\right) ^2ds.
 \end{equation}
 With the change of variable $s\mapsto n+1-s$, we get
 \[
B_n= \frac {c_H}n     \sum_{h=1}^n\int_{h}^{h+1} \left( \sum_{j=1}^{h}     \xi_j \Delta K_2( s-j)\right) ^2ds.
 \]
 
 Fix an integer $N\ge 1$ and let $M =[n/N  ]$ be the integer part of $n/N$. Then, $n\ge NM $.
 As a consequence,
\[
B_{n}  \ge  \frac {c_H}{n}     \sum_{i=0}^{N-1}    \sum_{h=iM+1}^{(i+1)M }  \int_{h}^{h+1} \left( \sum_{j=1}^{h}     \xi_j \Delta K_2(s-j)\right) ^2ds.
\]
Set
\[
B_{n}^i=\frac {c_H}{n}  \sum_{h=iM+1}^{(i+1)M }  \int_{h}^{h+1} \left( \sum_{j=1}^{h}     \xi_j \Delta K_2(s-j)\right) ^2ds,\quad i=0,\ldots,N-1.
\]
We are going to use the estimate
\begin{equation} \label{ineq}
(B_{n})^{-p} \le \prod_{i=0}^{N-1}  (B_{n}^i)^{-\frac pN}.
\end{equation}
Consider again the representation of the sequence $\xi$ as stochastic integrals  with respect to a Brownian motion $W$, and denote by $\{\mathcal{F}^W_t\}$ the filtration generated by the Brownian motion $W$. Then, 
\begin{eqnarray}
E[(B_{n})^{-p}] \le E\left[ \prod_{i=0}^{N-1}  (B_{n}^i)^{-\frac pN} \right]= E\left[    E[ (B_{n}^{N-1})^{-\frac pN}| \mathcal{F}^W_{(N-1)M}]\prod_{i=0}^{N-2}  (B_{n}^i)^{-\frac pN} \right].\label{induction}
\end{eqnarray}
Let us estimate the conditional expectation appearing in the right-hand side.
In the same spirit that (\ref{118bis}), 
it is immediate that
\[
E[ (B_{n}^{N-1})^{-\frac pN}| \mathcal{F}^W_{(N-1)M}]  \le  1 +\frac p N \int_0^1 P\left(  B_{n}^{N-1}  \leq x | \mathcal{F}^W_{(N-1)M} \right) x^{-\frac pN -1} dx.
\]
By  Carbery-Wright's inequality (\ref{ff1}) with $d=2$,
\begin{equation}\label{a}
P\left(  B_{n}^{N-1}  \leq x | \mathcal{F}^W_{(N-1)M} \right) \le   c \sqrt{x} \left[ E\left(   B_{n}^{N-1} | \mathcal{F}^W_{(N-1)M} \right) \right]^{-\frac1d}.
\end{equation}
The  conditional expectation  $E\left(   B_{n}^{N-1} | \mathcal{F}^W_{(N-1)M} \right)$ is given by
\begin{eqnarray}
\notag
E\left(   B_{n}^{N-1} | \mathcal{F}^W_{(N-1)M} \right)
=\frac {c_H}{n}   \sum_{h=(N-1)M+1}^{NM }  \int_{h}^{h+1} 
E\left[  \left( \sum_{j=1}^{h}     \xi_j \Delta K_2(s-j)\right) ^2 | \mathcal{F}^W_{(N-1)M} \right]ds.\\
\label{aa}
\end{eqnarray}
Taking into account that the sequence $\xi_j$ is Gaussian, the conditional expectation appearing in the above equation can be bounded below by the conditional variance which is {\it not} random. More precisely,
\begin{eqnarray}
&&E\left[  \left( \sum_{j=1}^{h}     \xi_j \Delta K_2(s-j)\right) ^2 | \mathcal{F}^W_{(N-1)M} \right]\ge  {\rm Var} \left[   \sum_{j=1}^{h}     \xi_j \Delta K_2(s-j)  \bigg| \mathcal{F}^W_{(N-1)M} \right]\notag\\
\qquad &&
=E\left[   \left( \sum_{j=(N-1)M+1}^{h}     \left(\int_{(N-1)M}^{j} \Delta K_1 (j, u) dW_u \right) \Delta K_2(s-j)\right) ^2    \right].\label{aaa}
\end{eqnarray}
By plugging (\ref{a}), (\ref{aa}) and (\ref{aaa}) into (\ref{induction}) and then by proceeding by induction with the other terms, we see that (\ref{claim}) will follow as soon as, for any $i=0,\ldots,N-1$,
\begin{equation} \label{eq4}
 \liminf_{M\rightarrow \infty}    \frac 1M \sum_{h=iM+1}^{(i+1)M}   \int_{h}^{h+1}  
E\left[   \left( \sum_{j=iM+1}^{h}     \left(\int_{iM}^{j} \Delta K_1(j,u) dW_u  \right) \Delta K_2(s-j)\right) ^2    \right]ds>0.
\end{equation}

 \medskip
 \noindent
 \textit{Proof of (\ref{eq4}).}   First we compute the expectation in (\ref{eq4}):
\begin{eqnarray*}
&& E\left[   \left( \sum_{j=iM+1}^{h}     \left(\int_{iM}^{j} \Delta K_1(j,u) dW_u  \right) \Delta K_2(s-j)\right) ^2    \right] \\
&& =\sum_{j,k=iM+1}^{h}  \Delta K_2(s-j)\Delta K_2(s-k) \int_{iM}^{j\wedge k}  \Delta K_1(j,u) \Delta K_1(k,u) du.
\end{eqnarray*}
Set
\[
\beta_{j,k}^M= \int_{iM}^{j\wedge k}  \Delta K_1(j,u) \Delta K_1(k,u) du.
\]
We can write, exchanging the order of the summation
\begin{eqnarray*}
&& \sum_{h=iM+1}^{(i+1)M}  \sum_{j,k=iM+1}^{h} \left(  \int_{h}^{h+1} \Delta K_2(s-j)\Delta K_2(s-k)ds  \right) \beta_{j,k}^M     \\
&&=\sum_{j,k=iM+1}^{(i+1)M }\sum_{h= j\vee k}^{(i+1)M} \left(  \int_{h}^{h+1} \Delta K_2(s-j)\Delta K_2(s-k)ds  \right) \beta_{j,k}^M=\sum_{j,k=iM+1}^{(i+1)M }  \alpha^M_{j,k}\, \beta_{j,k}^M,
\end{eqnarray*}
where
 \[
 \alpha^M_{j,k}=\int_ {j\vee k }^{(i+1)M+1} \Delta K_2(s-j)\Delta K_2(s-k)ds.
 \]
 Then we are interested in the liminf, as   $M$  tends to infinity, of
 \[
\frac 1M  \sum_{j,k=iM+1}^{(i+1)M  }\alpha^M_{j,k}\beta_{j,k}^M.
\]
 We make now the change of indices $j\mapsto j-iM$ and $k\mapsto k-iM$ and we obtain the expression
  \[
\Psi_M:= \frac 1M  \sum_{j,k= 1}^{M  }\widetilde{\alpha}^M_{j,k}\widetilde{\beta}_{j,k}^M,
\]
where
\begin{eqnarray*}
\widetilde{\beta}_{j,k}^M&=& \int_{iM}^{(j+iM)\wedge (k+iM)}  \Delta K_1(j+iM,u) \Delta K_1(k+iM,u) du\\
&=& \int_{0}^{j\wedge k}  \Delta K_1(j+iM,u+iM) \Delta K_1(k+iM,u+iM) du,
\end{eqnarray*}
and
  \[
\widetilde{ \alpha}^M_{j,k}=\int_ {j\vee k }^{ M+1} \Delta K_2(s-j)\Delta K_2(s-k)ds.
 \]
 
 \medskip
 {\it Step 1:} Case $H>\frac 12$. In this case $\Delta K_1$ and $\Delta K_2$ are nonnegative. On the other hand, by \cite[(5.10)]{nualartbook}, we have
 \[
 \frac {\partial K_1}{\partial t} (t,s)= c_H (t/s)^{H-\frac 12} (t-s) ^{H-\frac 32},
 \]
 where $c_H=\sqrt{\frac{H(2H-1)}{\beta(2-2H,H-1/2)}}$.
 Therefore, assuming $j\ge k\ge 2$,
\begin{eqnarray*}
 \widetilde{\beta}_{j,k}^M  &\ge&  \int_{1}^{  k-1}  \Delta K_1(j+iM,u+iM) \Delta K_1(k+iM,u+iM) du \\
& \ge&  c_H^2  \int_{1}^{  k-1}   \left( \int_{j+iM-1} ^{j+iM} \left( \frac x{u+iM} \right) ^{H-\frac 12} (x-u-iM) ^{H-\frac 32} dx \right)  \\
&&\quad \quad \quad \times   \left( \int_{k+iM-1} ^{k+iM} \left( \frac y{u+iM} \right) ^{H-\frac 12} (y-u-iM) ^{H-\frac 32} dy \right) du\\
&\ge & c_H^2  \int_{1}^{  k-1}   \left( \frac {(j+iM-1)(k+iM-1)}{(u+iM)^2} \right) ^{H-\frac 12}  (j-u) ^{H-\frac 32} (k-u) ^{H-\frac 32}du.
\end{eqnarray*}
 The term  $\frac{(j+iM-1)(k+iM-1)}{(u+iM)^2}$ is lower bounded by 1. 
  Therefore
\begin{eqnarray*}
  \widetilde{\beta}_{j,k}^M  & \ge&    c_H^2  \int_{1}^{  k-1}(j-u) ^{H-\frac 32} (k-u) ^{H-\frac 32}du
  =  c_H^2  \int_{1}^{  k-1 }(x+j-k ) ^{H-\frac 32} x ^{H-\frac 32}dx \\
  &\ge& c_H^2  \int_{1}^{  k-1 }(x+j-k ) ^{2H-3}  dx = C \left((j-k +1) ^{2H-2}- (j-1)^{2H-2} \right).
\end{eqnarray*}
By similar arguments we obtain, assuming again $j\ge k\ge 2$,
\begin{eqnarray*}
 \widetilde{\alpha}^M_{j,k}&\ge & \int_ {j +2 }^{ M+1} \Delta K_2(s-j)\Delta K_2(s-k)ds\\
 &\ge&  (H-\frac 12)  \int_ {j +2 }^{ M+1}  (s-j)^{H-\frac 32} (s-k) ^{H-\frac 32} ds
=   (H-\frac 12)  \int_ { 2 }^{ M+1-j}  x^{H-\frac 32}  (x+j-k) ^{H-\frac 32} dx\\
 &\ge & (H-\frac 12)  \int_ { 2 }^{ M+1-j}    (x+j-k) ^{2H-3} dx
 = C \left(  (j-k +2)^{2H-2} -(  M+1-k)^{2H-2})\right).
  \end{eqnarray*}
 Therefore,
\begin{eqnarray*}
 \Psi_M&\ge&  \frac 1M \sum_{2\le k\le j \le M}  \widetilde{ \alpha}^M_{j,k}  \widetilde{\beta}_{j,k}^M \\
 &\ge&  \frac C M \sum_{k=2}^M\sum_{\ell=0}^{M-k} \left((\ell +1) ^{2H-2}- (\ell+k-1)^{2H-2} \right)\left(  (\ell +2)^{2H-2} -(  M+1-k)^{2H-2})\right).
  \end{eqnarray*}
 For the first  term we obtain
 \[
 \frac 1M \sum_{k=2}^M\sum_{\ell=0}^{M-k} (\ell +1) ^{2H-2}(\ell +2) ^{2H-2} \ge
 \frac 1M \sum_{k=2}^M\sum_{\ell=0}^{M-k} (\ell +2) ^{4H-4} \ge \frac{M-1}M 2^{4H-4},
 \]
 which converges to a positive constant as $M$ tends to infinity. It is easy to check that the other terms in the above expression converge to zero. For instance,
 \begin{eqnarray*}
  \frac 1M  \sum_{k=2}^M\sum_{\ell=0}^{M-k}  (\ell +2) ^{2H-2}(\ell +k-1) ^{2H-2}
 & =& \frac 1M  \sum_{\ell=0}^{M-2}  (\ell +2) ^{2H-2} \sum_{k=2}^{M-\ell} (\ell +k-1) ^{2H-2}\\
 &\leq&\frac1M\left(\sum_{\ell=1}^{M}\ell^{2H-2}\right)^2
 \end{eqnarray*}
and this last quantity behaves as $M^{4H-3}$, which converges to zero because $H<\frac 34$. A similar analysis can be done for the other terms.
 
 \medskip
  {\it Step 2:} Case $H<\frac 12$.  In this case, see \cite[(5.23)]{nualartbook}, we have that
  \[
 \frac {\partial K_1}{\partial t} (t,s)= c_H (H-\frac 12) (t/s)^{H-\frac 12} (t-s) ^{H-\frac 32} \quad\mbox{(with $c_H=\sqrt{\frac{2H}{(1-2H)\beta(1-2H,H+1/2)}}$)}
 \]
 is negative. Therefore $\Delta K_1(j,s)$ is negative if $s<j-1$ and positive if $j-1\leq s\leq j$. Also,
   \[
 \frac {\partial K_2}{\partial t} (u)=  (H-\frac 12) u ^{H-\frac 32}
 \]
and $\Delta K_2(u)$ is negative if $u\geq 1$ and positive if $u<1$. Then, it suffices to show that the negative terms do not contribute to the limit, and once we get rid of these negative terms, we can get a lower bound as  in the case $H>\frac 12$. When $j=k$, the integrands in the definition of $\widetilde{\beta}_{j,k}^M$ and  $\widetilde{\alpha}_{j,k}^M$ are nonnegative. On the other hand, for $j\ge k+1$, we can write
\[
\widetilde{\alpha}_{j,k}^M=\widetilde{\alpha}_{j,k,1}^M+\widetilde{\alpha}_{j,k,2}^M,
\]
where
\[
\widetilde{\alpha}_{j,k,1}^M=\int_j^{j+1} (s-j)^{H-\frac 12} \left[ (s-k)^{H-\frac 12}-(s-k-1)^{H-\frac 12} \right]ds\le 0,
\]
and
\[
\widetilde{\alpha}_{j,k,2}^M=\int_{j+1}^{M+1}  \Delta K_2(s-j) \Delta K_2(s-k) ds\ge 0.
\]
Similarly
\[
\widetilde{\beta}_{j,k}^M=\widetilde{\beta}_{j,k,1}^M+\widetilde{\beta}_{j,k,2}^M,
\]
where
\[
\widetilde{\beta}_{j,k,1}^M=\int_{k-1}^{k}  K_1(k+iM,u+iM) \Delta K_1(j+iM,u+iM)du\le 0,
\]
and
\[
\widetilde{\beta}_{j,k,2}^M=\int_{0}^{k-1} \Delta K_1(k+iM,u+iM) \Delta K_1(j+iM,u+iM)du\ge 0.
\]
In this way we obtain the decomposition
\begin{eqnarray*}
\Psi_M &=&\frac 1M \sum_{j=1}^M  \widetilde{\alpha}_{j,j}^M\widetilde{\beta}_{j,j}^M
+ \frac 2M \sum_{1\le k \le j-1 \le M-1}\left( \widetilde{\alpha}_{j,k,1}^M+\widetilde{\alpha}_{j,k,2}^M\right)
\left( \widetilde{\beta}_{j,k,1}^M+\widetilde{\beta}_{j,k,2}^M\right)\\
&=& \Psi^1_M + \Psi^2_M + \Psi_M^3,
\end{eqnarray*}
where
\[
\Psi^1_M=\frac 1M \sum_{j=1}^M  \widetilde{\alpha}_{j,j}^M\widetilde{\beta}_{j,j}^M
+ \frac 2M \sum_{1\le k \le j-1 \le M-1}\left( \widetilde{\alpha}_{j,k,1}^M\widetilde{\beta}_{j,k,1}^M+\widetilde{\alpha}_{j,k,2}^M\widetilde{\beta}_{j,k,2}^M \right),
\]
\[
\Psi^2_M= \frac 2M \sum_{1\le k \le j-1 \le M-1}  \widetilde{\alpha}_{j,k,1}^M\widetilde{\beta}_{j,k,2}^M,
\]
and
\[
\Psi^3_M= \frac 2M \sum_{1\le k \le j-1 \le M-1}  \widetilde{\alpha}_{j,k,2}^M\widetilde{\beta}_{j,k,1}^M.
\]
The term $\Psi^1_M$ is nonnegative and it can be bounded below as follows
\begin{eqnarray*}
\Psi^1_M &\ge & \frac 2M \sum_{2\le k\le j-1 \le M-1} \left( \int_{j+2}^{M+1} \Delta K_2(s-j) \Delta K_2(s-k)ds \right)\\
&& \times 
\left( \int_1^{k-1} \Delta K_1 (j+iM, u+iM) \Delta K_1(k+iM, u+iM) du \right).
\end{eqnarray*}
By the same arguments as in the case $H>\frac 12$ we can show that $\liminf_{M\rightarrow \infty} \Psi^1_M >0$.
Therefore, it suffices to show that
\begin{equation} \label{c1}
\lim_{M\rightarrow \infty} 
\Psi^2_M=0,
\end{equation}
and
\begin{equation} \label{c2}
\lim_{M\rightarrow \infty} 
\Psi^3_M=0.
\end{equation}
These limits are based on the following estimates. One one hand,  $\widetilde{\beta}_{j,k,2}^M$ and  $\widetilde{\alpha}_{j,k,2}^M$ are uniformly bounded:
\begin{eqnarray*}
\widetilde{\beta}_{j,k,2}^M&=&\int_{iM}^{iM+k-1} \Delta K_1(k+iM,u) \Delta K_1(j+iM,u)du \\
&\le& \sqrt{\int_{iM}^{iM+k-1}  \Delta K_1(k+iM,u)  ^2du} \, \sqrt{\int_{iM}^{iM+k-1}  \Delta K_1(j+iM,u)  ^2du}\\
&\le& \sqrt{\int_{0}^{k+iM}  \Delta K_1(k+iM,u)  ^2du}\,  \sqrt{\int_{0}^{j+iM}  \Delta K_1(j+iM,u)  ^2du}\\
&=& \sqrt{E( |B_{k+iM}-B_{k+iM-1}|^2)} \sqrt{E( |B_{j+iM}-B_{j+iM-1}|^2)} =1,
\end{eqnarray*}
and
\begin{eqnarray*}
\widetilde{\alpha}_{j,k,2}^M
&=& \int_{j+1}^{M+1} [ (s-j)^{H-\frac 12} -(s-j-1)^{H-\frac 12}]  [ (s-k)^{H-\frac 12} -(s-k-1)^{H-\frac 12}]ds\\
&=& \int_1^{M+1-j} [x^{H-\frac 12} -(x-1) ^{H-\frac 12}] [(x+j-k)^{H-\frac 12} -(x+j-k-1) ^{H-\frac 12}]dx \\
&\le & C+  \int_2^{M+1-j} [x^{H-\frac 12} -(x-1) ^{H-\frac 12}] [(x+j-k)^{H-\frac 12} -(x+j-k-1) ^{H-\frac 12}]dx \\
&\le & C+ C_H\int_2^{M+1-j} (x-1) ^{H-\frac 32}  (x+j-k-1) ^{H-\frac 32}dx\\
&\le & C+ C_H\int_2^{M+1-j} (x-1) ^{2H-3}  dx \le C'.
\end{eqnarray*}
On the other hand, we have
\begin{eqnarray*}
|\widetilde{\alpha}_{j,k,1}^M| &\le&  \int_0^1 x^{H-\frac 12} \left[ j-k +x-1)^{H-\frac 12} - (j-k+x) ^{H-\frac 12} \right] dx\\
&=&\frac{1}{\frac32-H}\int_0^1 x^{H-\frac12}\left(\int_{j-k-1}^{j-k}(u+x)^{H-\frac32}du\right)dx
\le \frac{1}{(\frac32-H)(\frac12+H)} (j-k-1)^{H-\frac 32},
\end{eqnarray*}
 and
\begin{eqnarray*}
|\widetilde{\beta}_{j,k,1}^M| &\le& \int_{k-1}^{k} K_1(k+iM,u+iM)| \Delta K_1(j+iM,u+iM)|du\\
&=& c_H \int_{k-1}^k K_1(k+iM, u+iM)  \left( \int_{j+iM-1}^{j+iM} \left( \frac x{u+iM} \right)^{H-\frac 12} (x-u-iM)^{H-\frac 32} dx \right) du \\
&\le &  c_H \left( \frac{j+iM-1}{ k+iM} \right) ^{H-\frac 12} \int_{k-1}^k K_1(k+iM, u+iM) \left( \int_{j-1}^{j}(x-u)^{H-\frac 32} dx  \right) du \\
&\le & C (j-k-1)^{H-\frac 32},
\end{eqnarray*}
because  $ \frac{j+iM-1}{ k+iM}\geq 1$, the integral $\int_{j-1}^{j}(x-u)^{H-\frac 32} dx $ is bounded by $(j-k-1)^{H-\frac 32}$ uniformly in $u\in [k-1,k]$, and
\[
\sup_{k,M} \int_{k-1}^k K_1(k+iM, u+iM) du <\infty, 
 \]
 as it can be easily checked from the expression of $K_1$.
 Finally,
 \[
 \lim_{M\rightarrow \infty}  \frac 1M \sum_{1\le k \le j-2 \le M-2} (j-k-1)^{H-\frac 32}=0,
 \]
 which implies both (\ref{c1}) and (\ref{c2}).

\end{document}